\documentclass{sl}     

\usepackage{slsec}     
\usepackage{stud_log} 
\usepackage{slfoot}
\usepackage{slthm}     
                         
\def\state#1#2{

\ 

\noindent{\bf #1} \emph{#2}

}

\usepackage{amsmath}
\usepackage{amssymb}

\usepackage{enumerate}

\def\kn{\kern.1em}

\renewcommand{\theequation}{\Roman{equation}}

\theoremstyle{definition}

\HeadingsInfo{Slava Meskhi}{Injectives in the variety generated by a finite subdirectly irreducible Heyting algebra with involution }

\begin{document}

\setcounter{page}{1}     

\hyphenation{Qua-cken-bush}
 

\AuthorTitle{Slava Meskhi}{Injectives in the variety generated by a finite subdirectly irreducible Heyting algebra with involution}







\begin{abstract}We prove that any finite subdirectly irreducible Heyting algebra with involution is quasi-primal, and that injective algebras in the variety generated by a finite subdirectly irreducible Heyting algebra are precisely diagonal subalgebras of some direct power of this algebra, which are complete as lattices.
\end{abstract}



\section{Introduction and overview} 

This note is devoted to Heyting algebras with involution, i.e. the automorphism of order two. Using the one-to-one correspondence between the set of all congruences of a Heyting algebra with involution and the set of all involutive filters of the algebra, we prove that any finite Heyting algebra with involution is quasi-primal. R.~W.~Quackenbush in [12] gives in terms of Boolean powers a nice characterization of injective algebras in a variety generated by a quasi-primal algebra.

The main purpose of this paper is to illustrate how the result of R.~W.~Quackenbush can be applied to Heyting algebras with involution. Combining the Quackenbush's result with our description of Boolean powers of quasi-primal algebras with lattice reducts, we prove that injectives in the variety generated by a finite Heyting algebra $\textbf{A}$ are diagonal subalgebras of some direct power of $\textbf{A}$, which are complete as lattices.

The sketch of the main points of this paper can be found in [10], which is devoted to some discriminator subclass of Heyting algebras with involution, namely, Heyting algebras with regular involution. This subclass of Heyting algebras with involution is semi-simple and contains simple algebras of arbitrary cardinality [10]. It is easy to see that if a variety has the minimal nontrivial algebra embeddable in any nontrivial algebra of this variety, and for any cardinal $\aleph$ there exists a simple algebra of cardinality $\aleph$, then there are no nontrivial injectives in this variety. That is why we focus our interest on the description of injectives in varieties generated by some finite Heyting algebras with involution.

Our motivation to study Heyting algebras with involution comes from many-valued logic. In [7] we showed that linearly ordered Heyting algebras with involution are the algebraic models of fuzzy propositional logic. Generalization of the methods of the paper mentioned lead us to the notion of Heyting algebra with regular involution [8]. It must be mentioned that Heyting algebras with involution are considered in [2] under the name "involutive Skolem algebras", and in [3] the duality theory for these algebras is announced.

\section{Preliminaries} %

In this section we will introduce the definitions and we will provide some equations which are satisfied in the variety of Heyting algebras with involution.

An algebra $\textbf{A} =\langle A, \lor, \land , \to , \neg , 0, 1 \rangle$ is a \emph{Heyting algebra} if  $\langle A, \lor, \land, 0, 1 \rangle$  is a bounded lattice; $\to$ is a binary operation on $A$ defined by $a\to  b =\sup\ \{x: a\land  x\le  b\}$, and $\neg$  is an unary operation of pseudocomplementation defined by $\neg  a = a\to  0$, or, in other terms, $\neg a =\sup\ \{x: a\land  x = 0\}$. We denote by  $\textbf{H}$ the variety of all Heyting algebras.

An algebra $\textbf{A} =\langle A, \lor, \land ,\to ,\neg ,\sim , 0, 1 \rangle$  is a \emph{Heyting algebra with involution} if $\langle A, \lor, \land, \to, \neg , 0, 1 \rangle$  is a Heyting algebra and $\sim$  is an unary operation for which the following equations are satisfied:

i1. $\sim  (a \lor b) =\sim  a\land\sim  b$

i2. $\sim\sim  a = a$

We denote by $\textbf{HI}$ the variety of all Heyting algebra with involution.
In other terms, an algebra $\textbf{A} =\langle A, \lor, \land, \to, \neg ,\sim , 0, 1 \rangle$ belongs to $\textbf{HI}$ if $\langle A, \lor, \land, \to, \neg , 0, 1 \rangle$ is a Heyting algebra and $\langle A, \lor, \land, \sim, 0, 1 \rangle$  is a De Morgan algebra. Recall that in any De Morgan algebra the following holds:

i3. $\sim  (a \land  b) =\sim  a \lor \sim  b$

i4. $\sim  0 = 1;\sim  1 = 0$

i5. $a \leq  b  \Rightarrow \sim  b \leq  \sim  a$

Using the properties of Heyting algebras and De Morgan algebras it is easy to prove the following

\state{Lemma 2.1.}{For any Heyting algebra with involution the following holds:}

i3. $\sim  (a \land  b) =\sim  a \lor \sim  b$

i4. $\sim  0 = 1;\sim  1 = 0$

i5. $a \leq  b \Rightarrow \sim  b \leq \sim  a$
\begin{proof}
i6 follows from i3 and the fact that in any Heyting algebra $\neg  (a \lor b) = \neg  a \land \neg  b$.

i7 follows from i4 and the fact that in any Heyting algebra $\neg  1 = 0$; $\neg  0 = 1$.

i8 follows from i5 and the fact that in any Heyting algebra $a \leq  b \Rightarrow \neg  b \leq \neg  a$.

Let us check i9.

$(a \to b) \land a \le b$ (the definition of $\to$ );

$\sim b \le \sim (a\to b) \lor \sim a$ (i5, i3);

$\neg \sim (a \to b) \land \sim b \le (\neg \sim (a \to b) \land \sim (a\to b)) \lor \neg \sim (a\to b) \land \sim a)$ (distributive lattice operations $\land$ and $\lor$ properties);

but $\neg \sim (a \to b) \land \sim (a \to b) = 0$ (the definition of $\neg$),

so $\neg \sim (a \to b) \land \sim b \le \neg \sim (a \to b) \land \sim a$ and
$\neg \sim (a \to b) \land \sim b \le \sim a$;
$\neg \sim (a \to b) \le \sim b \to \sim a$

\end{proof}

\section{Subdirectly irreducible algebras}

In this section we will give the characterization of a finite subdirectly irreducible Heyting algebra with involution in terms of involutive centre of this algebra. Basic universal algebraic definitions of this section can be found in [1], [5].

A subalgebra $\textbf{A}$ of the direct product of algebras $\{\textbf{A}_i: i\in I\}$ is said to be a \emph{subdirect product} of $\{\textbf{A}_i: i\in I\}$ if for any $i\in I$ the projection $\pi_i(A) = A_i$.

An algebra $\textbf{A}$ is said to be \emph{subdirectly irreducible} if for any set of congruences $\{\Theta_i: i\in I\}$ from $\bigcap\{\Theta_i: i\in I\} =\Delta$  it follows that $\Theta_i =\Delta$  for some $i\in I$, where $\Delta$  is the diagonal relation. In other words, $\textbf{A}$ is subdirectly irreducible iff $\textbf{A}$ either has only one congruence or there exists the smallest element in the set of all congruences not equal to $\Delta$.

For any lattice $A$, a subset $F$ of elements from $A$ is said to be a \emph{filter} if:

1) $a,b\in F$ $\Rightarrow$  $a\land  b\in F$;

2) $a\in F$ and $a\leq  b$ $\Rightarrow$ $b\in F$.

It is well known that there exists a one-to-one correspondence between the set of all congruences of a Heyting algebra  $\textbf{A}$ and the set of all filters of $\textbf{A}$ (see, for example [13]). The correspondence is defined as follows:
for a filter $F$ of $\textbf{A}$, the relation $\Theta (F)$, where $\Theta (F)(a,b)\iff  (a\to  b)\land  (b\to  a)\in F$, is a congruence of $\textbf{A}$, and for a congruence $\Theta$  of $\textbf{A}$ the set $F(\Theta ) =\{a\in A:\Theta (a,1)\}$ is a filter.

Note that the smallest congruence $\Delta$  corresponds to the unit filter $\{1\}$.

A filter $F$ of a Heyting algebra with involution $\textbf{A}$ is said to be \emph{involutive} if for any $a\in A$

$a\in F$ $\Rightarrow$  $\neg\sim  a\in F$

\state{Lemma 3.1.}{For $\textbf{A}\in\textbf{HI}$ there exists a one-to-one correspondence between the set of all congruences of $\textbf{A}$ and the set of all involutive filters of $\textbf{A}$.}
\begin{proof}
Let $F$ be an involutive filter of $\textbf{A}$. Then the relation

$$\Theta (F)(a,b)\ \iff\ (a\to  b)\land  (b\to  a)\in F$$

is an equivalence, which preserves all Heyting operations (see, for example [13]). We will show that $\Theta (F)$ preserves the involution. Let $\Theta (F)(a,b)$. Then

$(a\to  b)\in F\ \Rightarrow\ \neg\sim  (a\to  b)\in F\ \Rightarrow\  (\sim  b\to\sim  a)\in F;$

$(b\to  a)\in F\ \Rightarrow\  \neg\sim  (b\to  a)\in F\ \Rightarrow\  (\sim  a\to\sim  b)\in F$ 

using the property of an involutive filter and Lemma 2.1 (i9). Thus, we have $\Theta (F)(\sim  a,\sim  b)$.

Now, let $\Theta$ be a congruence of $\textbf{A}$. Then the set $F(\Theta ) =\{a\in A:\Theta (a,1)\}$ is a filter, and from Lemma 2.1 (i7) it follows that $F(\Theta )$ is an involutive filter.

\end{proof}

For any $\textbf{A}$ from $\textbf{HI}$ and a subset $B$ of $A$, the intersection of all involutive filters containing $B$ (for example, $A$) is the smallest involutive filter containing $B$. Such a filter is said to be generated by $B$, and is denoted by $F[B]$.

We will use the following abbreviation $[\neg\sim  n]a =\neg\sim(\neg\sim(... (\neg\sim a)...))$ where $n-1$ is the number of left brackets.

\state{Lemma 3.2.}{For $\textbf{A}$ from $\textbf{HI}$ and a subset $B$ of $A$, the following condition is satisfied:
$$
\begin{aligned}
&\textrm{$a\in F[B]$ iff}\\
&\textrm{there exist $a_1, a_2,..., a_m\in B$ and $n(a_1), n(a_2),..., n(a_m)$ such that}\\
&[\neg\sim  n(a_1)]a_1\land  [\neg\sim  n(a_2)]a_2\land...\land  [\neg\sim  n(a_m)]a_m\le a.
\end{aligned}
\leqno{(*)}
$$
}
\begin{proof}
It is evident that if for $a\in A$ the condition (*) holds and $a\leq  b$ then for $b$ the condition (*) holds. Moreover, from Lemma 2.1 (i6) follows that if for $a,b\in A$ the condition (*) holds, then it holds for $a\land  b$. Thus, the set of all elements, for which the condition (*) holds is a filter, and from Lemma 2.1 (i8) follows that this filter is an involutive one.
\end{proof}

Let $\textbf{A}\in\textbf{HI}$. The set of all elements $a\in A$ such that $\neg  a =\sim  a$ we call an involutive center of $\textbf{A}$ and denote it by $IC(A)$.

\state{Lemma 3.3.}{ For a Heyting algebra with involution $\langle A, \lor, \land ,\to ,\neg, \sim , 0, 1 \rangle$, the involutive center $\langle IC(A) , \lor, \land, \to, \neg , 0, 1 \rangle$  is a Boolean algebra.}
\begin{proof}
For any $a\in IC(A)$ 
$$a \lor \neg  a = a \lor \sim  a = \sim \sim  a \lor \sim  a = \sim  (\sim a \land  a) =\sim  (\neg a \land  a) = \sim  0 = 1$$
so $\neg$ is the Boolean complement operation
\\ Let us show that $IC(A)$ is closed under the operatins of $\lor$, $\land$ , $\neg$. Let $a,b\in IC(A)$. Then  $\sim  (a \lor b)$ is the complement of the element $a \lor b$.

Indeed,
\begin{align*}
(a \lor b)\land\sim  (a \lor b)& = (a \land \sim  (a \lor b)) \lor (b \land \sim  (a \lor b))\\
&= (a \land \sim a \land \sim  b) \lor (b \land \sim a \land \sim  b)\\
&= (a \land \neg a \land \sim  b) \lor (b \land \sim a \land \neg  b) = 0 \land  0 = 0.\\
\\
(a \lor b) \lor \sim (a \lor b)& = (a \lor b) \lor (\sim a \land \sim b)\\
&= (a \lor b \lor \sim a)\land  (a \lor b \lor \sim  b)\\
&= (a \lor b \lor \neg a) \land  (a \lor b \lor \neg  b) = 1 \lor 1 = 1.
\end{align*}
But the complement of an arbitrary element is the pseudocomplement of this element, and from the uniqueness of the pseudocomplement, we have:

$\neg  (a \lor b) =\sim  (a \lor b)$ and $\neg  (a \land  b) =\sim (a \land  b)$.

Let us show that for $a\in IC(A)$, $\neg  a\in IC(A)$.
$\neg\sim  a =\sim \sim  a = a$. At the same time, $\neg a$ is the complement of $a$, and therefore, the complement of $\neg a$ is $a$, i.e. $\neg\neg  a = a$. Thus $\neg \sim  a =\neg\neg  a$.
\end{proof}

\state{Theorem 3.4.}{For an arbitrary finite Heyting algebra with involution $\textbf{A} =\langle A, \lor,\land, \to,\neg, \sim, 0, 1 \rangle$  there exists a one-to-one correspondence between the set of all involutive filters of $A$ and the set of all filters of the Boolean algebra $\langle IC(A), \lor,\land, \to, \neg,  0, 1 \rangle$.}
\begin{proof}
Let $F$ be an involutive filter of $\textbf{A}$. Then the intersection of $F$ and $IC(A)$ is not empty ($1\in IC(A)$) and is a filter of $IC(A)$. Conversely, for a filter of the Boolean algebra $IC(A)$, the involutive filter of $\textbf{A}$ is constructed as it is shown in Lemma 3.2. All that we have to show is that two different involutive filters cannot have the same intersection with $IC(A)$. Suppose that such filters exist, i.~e. $F$ and $G$ are involutive filters of $\textbf{A}$, $F\cap IC(A) = G\cap IC(A)$ and for some $a\in A$, $a\in F$ and $a\notin G$.  Then $\neg\sim  a\in F$. We consider two cases:
\begin{enumerate}
\item
$a \leq \neg \sim a$. Therefore $a \land \sim a \leq \neg \sim a \land \sim  a = 0$. Then $a \lor \sim  a = 1$. So $\sim a = \neg  a$ and $a\in IC(A)$, which contradicts $F\cap IC(A) = G\cap IC(A)$.
\item
It is not true that $a \leq \neg \sim a$. Then $a \land \neg \sim a < a$. In this case there exists a new element ($a \land \neg \sim a$) which is strictly less than $a$, and ($a \land \neg \sim a)\in F$ and $(a\land\neg\sim  a)\notin G$. For this element we can consider the same two cases. As the result, we either will get, on some finite step, the contradiction to case (1), or will get the descending chain of elements, which, due to the fact that A is a finite algebra, will include 0. But $0\in IC(A)$. Then $F = G = A$.
\end{enumerate}
This contradiction proves the theorem.
\end{proof}
It is a well known fact that there is only one nontrivial subdirectly irreducible Boolean algebra ${\mathbf2}=\{0,1\}$ (see for example [4]). This yields, together with Lemma 3.1, Lemma 3.3, and Theorem 3.4, the following:

\state{Corollary 3.5.}{For a nontrivial finite algebra $\textbf{A}\in \textbf{HI}$ the following conditions are equivalent:
\begin{itemize}
\item[1.]
$\textbf{A}$ is subdirectly irreducible;
\item[2.]
$IC(A)$ is subdirectly irreducible;
\item[3.]
There is no $a\in A$ such that $a\ne 1$, $a\ne 0$, $\neg  a =\sim  a$
\end{itemize}
}

\section{The discriminator and the killer}

In this section we will show that any finite subdirectly irreducible Heyting algebra with involution is quasi-primal. Thereto we will prove that for Heyting algebras with involution the notion of discriminator is equivalent to the notion of killer, which is introduced in this section.

Recall some definitions. All of them except the definition of killer can be found in [5] and in [1].

For an algebra $\textbf{A}$ the function $t:A^3\to A$ is said to be a \emph{discriminator function} on $A$ if for any $a,b,c\in A$
$$
t(a,b,c) =
\begin{cases}
a&\textrm{if }a\ne b\\
c&\textrm{if }a = b
\end{cases}
$$

A ternary term $t(x,y,z)$, which represents the discriminator function on $\textbf{A}$, is said to be a \emph{discrimanator} for $\textbf{A}$.
A variety $\textbf{V}$ is said to be a \emph{discriminator variety} if there exists a class of algebras $\textbf{K}\subseteq\textbf{V}$ such that $\textbf{V}=V(\textbf{K})$ (i.~e. $\textbf{V}$ is generated by $\textbf{K}$) and there exists a common ternary discriminator for algebras from $\textbf{K}$.  Examples of discriminator varieties are: Boolean algebras, cylindric algebras of dimension $n$, Post algebras of order $n$, \L{}ukasiewicz algebras of order $n$. We will show that varieties generated by finite Heyting algebras with involution have the same property.

For an algebra $\textbf{A}\in \textbf{HI}$ a function $k:A\to A$ is said to be a \emph{killer} on $\textbf{A}$ if for any $a\in A$
$$
	k(a) =
\begin{cases}
0&\textrm{if $a\ne 1$}\\
1&\textrm{if $a = 1$}
\end{cases}
$$
 An unary term $k(x)$, which represents the killer function on $\textbf{A}$, is said to be a \emph{killer} for $\textbf{A}$.

\state{Lemma 4.1.}{Let $\textbf{K}$ be any finite set of finite subdirectly irreducible Heyting algebras with involution. Then there exists a common killer for algebras from $\textbf{A}$.}

\begin{proof}
Let N be the maximal length of all maximal chains of algebras from $\textbf{K}$. By $k_1(x)$ we denote $x\land\neg\sim  x$, and for any $i\leq  N$ let $k_i(x) = k_{i-1}(x)\land\neg\sim  k_{i-1}(x)$. We will show that $k_N(x)$ is the common killer for algebras from $\textbf{K}$.

Let  $\textbf{A}\in\textbf{K}$ and $a\in A$ with $a\ne 0$. It is clear that $k_1(a) = a \land \neg \sim a \leq a$.

Suppose $a = a \land \neg \sim a$. Then  $\sim a \land a = \sim a \land a\land \neg \sim a = 0 \land  a = 0$. Hence $\sim (\sim a \land a) =\sim  0 = 1$, therefore $a \lor \sim a = 1$. Then $\sim a = \neg a$, which is impossible for a nonzero element in a subdirectly irreducible algebra $\textbf{A}$ (Corollary 2.5.). So $k_1(a)$ is strictly below $a$. If $k_1(a) = 0$, then $k_i(a) = 0$ for any $i\leq N$. Otherwise, the same arguments proves that $k_i(a)$ is strictly less than $k_{i-1}(a)$, on some step $i$, from 1 to $N$, $k_i(a) = 0$.

It is evident that $k_i(1) = 1$ for any $i\leq  N$.
\end{proof}

\state{Lemma 4.2.}{Let $\textbf{K}$ be some class of Heyting algebras with involution. Then the following conditions are equivalent:
\begin{enumerate}
\item
there exists a common killer for algebras from $\textbf{K}$;
\item
there exists a common discrimanator for algebras from $\textbf{K}$.
\end{enumerate}
}
\begin{proof}
Let $k$ be a killer term for algebras from $\textbf{K}$. We will show that the following term
$$
t(x,y,z) = (k((x \to y) \land (y \to x)) \land z) \lor (\sim k((x \to y) \land  (y \to x)) \land  x)
$$
represents a discrimanator function for any algebra from $\textbf{K}$.

Let $a,b,c\in A$, for some algebra $\textbf{A}$ from $\textbf{K}$. First we consider the case $a\ne b$. In any Heyting algebra $(a \to b) \land  (b \to a) = 1$ iff $a = b$. Hence $(a \to b) \land  (b \to a) \ne 1$, $k((a \to b) \land (b \to a)) = 0$ and $\sim k((a \to b) \land (b\to a)) = 1$. Therefore $t(a,b,c) = 0 \lor a = a$.

Now we consider the case $a = b$. Then $(a \to b) \land (b\to a) = 1$. Hence $k((a \to b) \land  (b \to a)) = 1$, and $\sim k((a \to b) \land  (b \to a)) = 0$.

Hence, $t(a,b,c) = c \lor 0 = c$.

Thus the polynomial $t(x,y,z)$ discriminates y and represents a discriminator for algebras from $\textbf{K}$.

Now let $t$ be a discriminator polynomial for algebras from $\textbf{K}$. We will show that the polynomial $\sim t(1,x,0)$ is a killer for algebras from $\textbf{K}$. Let $a \in A$, for some algebra $\textbf{A}$ from $\textbf{K}$. If $a \ne 1$ then $t(1,a,0) = 1$, therefore $\sim  t(1,a,0) = 0$. If $a = 1$ then $t(1,a,0) = 0$ and therefore $\sim t(1,a,0) = 1$. This proves the Lemma.
\end{proof}

\state{Theorem 4.3.}{Let $\textbf{V}=V(\textbf{K})$  be a variety generated by some finite subset $\textbf{K}$ of finite subdirectly irreducible Heyting algebras with involution. Then $\textbf{V}$ is a discriminator variety.}

\begin{proof}
 By Lemma 4.1 there exists a common killer for algebras from $\textbf{K}$, and by Lemma 4.2 there exists a common discrimanator for algebras from $\textbf{K}$.
\end{proof}

Quasi-primal algebra, introduced by A. F. Pixley in [11] under the name "simple algebraic algebra", can be defined in terms of discrimanator. A finite nontrivial algebra $\textbf{A}$ is quasi-primal if there exists a discrimanator for $\textbf{A}$ (see [5]).

\state{Corollary 4.4.}{Any finite subdirectly irreducible Heyting algebra with involution is quasi-primal.}

\section{Injective algebras}
In this section we will show that injective algebras in the variety generated by a finite Heyting algebra with involution  $\textbf{A}$ are the diagonal subalgebras of some direct power of  $\textbf{A}$ which are complete as lattices.

An algebra $\textbf{A}$ from a class of algebras $\textbf{K}$ is said to be injective ($\textbf{K}$-injective), if for any pair $\textbf{B},\textbf{C} \in \textbf{K}$ such that $\textbf{B}$ is a subalgebra of $\textbf{C}$, any homomorphism $h:\textbf{B} \to \textbf{A}$ is extendable to some homomorphism $f: \textbf{C} \to \textbf{A}$ (i.e. the restriction of $f$ to $\textbf{B}$ coincides with $h$).

In [12] R.~W.~Quackenbush gives the following  nice characterisation of injective algebras in varieties generated by quasi-primal algebras.

\state{Theorem 5.1.}{Algebra $\textbf{A}$ is injective in the variety generated by some quasi-primal algebra $\textbf{C}$ iff $\textbf{A}$ is isomorphic to a Boolean power $\textbf{C}[\textbf{B}]$ for some complete Boolean algebra $\textbf{B}$.}
$ $

On the other hand, in [9] we proved that for quasi-primal algebras which have lattice operations in the signature, the notion of a complete Boolean power is reducible to the notion of a diagonal subalgebra of the direct power which is complete as a lattice:

\state{\bf Theorem 5.2.}{Let $\textbf{A}$ be a nontrivial quasi-primal algebra with the lattice operations in the signature. Then the following conditions are equivalent:
\begin{itemize}
\item[(i)]
an algebra $\textbf{C}$ is isomorphic to a Boolean power $\textbf{A}[\textbf{B}]$ for some complete Boolean algebra $\textbf{B}$.
\item[(ii)]
$\textbf{C}$ is complete as a lattice and $\textbf{C}$ is a diagonal subalgebra of a direct power $\textbf{A}^I$ for some $I$.
\end{itemize}
}
$ $
Combining Corollary 3.3 with Theorems 5.1 and 5.2, we get the following characterization of injective algebras in varieties generated by some finite Heyting algebra with involution:

\state{Theorem 5.3.}{Let $\textbf{A}$ be a nontrivial finite subdirectly irreducible Heyting algebra with involution and let $V(\textbf{A})$ be the variety generated by $\textbf{A}$. Then the following conditions on an algebra $\textbf{C}$ are equivalent:
\begin{itemize}
\item[(i)]
$\textbf{C}$ is injective in $V(\textbf{A})$.
\item[(ii)]
$\textbf{C}$ is complete as a lattice and $\textbf{C}$ is a diagonal subalgebra of a direct power $\textbf{A}^I$ for some $I$.
\end{itemize}
}
$ $
Note that if $\textbf{A}$ has only two elements, and both of them are constants, any subalgebra of the direct power $\textbf{A}^I$ is diagonal, and we get the well-known characterization of injective Boolean algebras (see for example [6], [4]).



\AuthorAdressEmail{Slava Meskhi}{Institute of Cybernetics\\
of the Georgian Technical University\\
Address\\Sandro Euli str. 5,
Tbilisi, 0186, Georgia}{slava.meskhil@gmail.com}

\end{document}